\documentclass[11pt,reqno]{amsart}

\usepackage[utf8]{inputenc}
\usepackage[T1]{fontenc}
\usepackage{lmodern}
\usepackage{microtype}
\usepackage[UKenglish]{babel}
\usepackage{amsmath,amsthm,amssymb,amsfonts}
\usepackage{mathrsfs}
\usepackage{mathtools}
\usepackage{enumitem}
\usepackage{xcolor}

\definecolor{dblue}{rgb}{0,0,.5}
\usepackage[pdfauthor={Matthias Paulsen},pdftitle={On the degree of algebraic cycles on hypersurfaces},colorlinks=true,linkcolor=dblue,citecolor=dblue,filecolor=dblue,menucolor=dblue,urlcolor=dblue]{hyperref}

\setlist[enumerate]{label=(\roman*)}

\newtheorem{theorem}{Theorem}
\newtheorem{conjecture}[theorem]{Conjecture}
\newtheorem{lemma}[theorem]{Lemma}
\newtheorem{proposition}[theorem]{Proposition}
\newtheorem{corollary}[theorem]{Corollary}
\theoremstyle{definition}
\newtheorem{remark}[theorem]{Remark}

\newcommand*{\PP}{\ensuremath{\mathbb{P}}}
\newcommand*{\CC}{\ensuremath{\mathbb{C}}}

\newcommand*{\QQ}{\ensuremath{\mathbb{Q}}}
\newcommand*{\ZZ}{\ensuremath{\mathbb{Z}}}
\newcommand*{\FF}{\ensuremath{\mathbb{F}}}
\newcommand*{\Ocal}{\ensuremath{\mathcal{O}}}
\newcommand*{\pr}{\ensuremath{\mathrm{pr}}}
\newcommand*{\dd}{\ensuremath{\mathrm{d}}}

\begin{document}

\title[On the degree of algebraic cycles on hypersurfaces]{On the degree of \\ algebraic cycles on hypersurfaces}
\author{Matthias Paulsen}
\address{Institute of Algebraic Geometry \\ Gottfried Wilhelm Leibniz Universität Hannover \\ Welfengarten~1 \\ 30167 Hannover \\ Germany}
\email{paulsen@math.uni-hannover.de}
\date{\today}

\begin{abstract}
Let $X\subset\PP^4$ be a very general hypersurface of degree $d\ge6$.
Griffiths and Harris conjectured in 1985 that the degree of every curve $C\subset X$ is divisible by $d$.
Despite substantial progress by Koll\'ar in 1991, this conjecture is not known for a single value of $d$.
Building on Koll\'ar's method, we prove this conjecture for infinitely many $d$, the smallest one being $d=5005$.
The set of these degrees~$d$ has positive density.
We also prove a higher-dimensional analogue of this result and construct smooth hypersurfaces defined over $\QQ$ that satisfy the conjecture.
\end{abstract}

\maketitle

\section{Introduction}

Throughout this paper, we work over the field $\CC$ of complex numbers.

In their famous work \cite{griffiths-harris}, Griffiths and Harris
made five conjectures about curves on a very general hypersurface $X\subset\PP^4$.
The weakest one is:
\begin{conjecture}[Griffiths--Harris]\label{conj:gh}
Let $X\subset\PP^4$ be a very general hypersurface of degree $d\ge6$.
Then the degree of every curve $C\subset X$ is divisible by $d$.
\end{conjecture}
This would follow from the stronger conjectures that $C$ is algebraically or even rationally equivalent to a multiple of a plane section. Their strongest conjecture, stating that $C$ is a complete intersection with a surface in $\PP^4$, was disproven by Voisin \cite{voisin}.
In contrast, Wu proved in \cite{wu} that every curve $C\subset X$ of degree at most $2d-2$ is a complete intersection with a surface in $\PP^4$.
In particular, the degree of every curve $C\subset X$ is at least $d$.
Nevertheless, Conjecture~\ref{conj:gh} is still open in any degree~$d$.

More generally, one might conjecture:
\begin{conjecture}\label{conj:higher}
Let $n\ge1$ be an integer, and let $X\subset\PP^{n+1}$ be a very general hypersurface of degree $d\ge2n$. Then the degree of every positive-dimensional closed subvariety $Z\subset X$ is divisible by $d$.
\end{conjecture}
Note that this conjecture is wrong if $1<d<2n$ because a hypersurface $X\subset\PP^{n+1}$ of degree $d<2n$ contains a line.
The case $n=1$ of this conjecture is trivial, and the case $n=2$ follows from the Noether--Lefschetz theorem \cite{lefschetz}. For $n\ge3$, however, Conjecture~\ref{conj:higher} is not known for a single $d$ yet.

\subsection{Koll\'ar's method}
For a very general hypersurface $X\subset\PP^4$ of degree~$d$, let us write
\[ f_3(d) = \gcd\{\deg C\mid C\subset X \text{ curve} \} \;. \]
Conjecture~\ref{conj:gh} states that $f_3(d)=d$ for all $d\ge6$.

The Trento examples \cite{trento-examples}, mostly due to Koll\'ar, achieve substantial progress towards Conjecture~\ref{conj:gh}
via specialization arguments.
By degenerating a very general hypersurface $X\subset\PP^4$ into a singular projection of a smooth projective threefold $Y$,
the following results are obtained:
\begin{enumerate}
\item\label{it:ddd} $d\mid 6\cdot f_3(d^3)$ for all $d\ge1$ (Koll\'ar, see also \cite[section~2]{soule-voisin})
\item\label{it:3dd} $d\mid 6\cdot f_3(3d^2)$ for all $d\ge4$ (Koll\'ar)
\item\label{it:6d} $d\mid 2\cdot f_3(6d)$ for all $d\ge9$ (van Geemen, improved by \cite{debarre-hulek-spandaw})
\end{enumerate}
These results naturally generalize to arbitrary dimension $n\ge3$ (see section~\ref{sect:trento}),
but they do not prove Conjecture~\ref{conj:higher} for any $d\ge2n$.

\subsection{Main results}
The main purpose of this article is to show the following:
\begin{theorem}\label{thm:main}
Let $n\ge3$ be an integer. Then there exists a set of degrees~$d$ with positive density such that Conjecture~\ref{conj:higher} is true in degree~$d$.
\end{theorem}
In particular, the case $n=3$ proves the conjecture of Griffiths and Harris for infinitely many degrees~$d$.
The smallest of them is $d=5005$.
Hence, $d=5005$ is the first degree where Conjecture~\ref{conj:gh} is currently known.

For a projective variety $X$, let us introduce the group
\[ Z^{2c}(X) = \frac{H^{c,c}(X,\ZZ)}{\left<\text{alg.\ classes}\right>} \;, \]
which measures the failure of the integral Hodge conjecture on $X$ in codimension~$c$.
As a consequence of Theorem~\ref{thm:main}, we get:
\begin{corollary}\label{cor:Z}
Let $n\ge3$ be an integer.
Then there exists a set of degrees~$d$ with positive density such that
\[ Z^{2c}(X)\cong\ZZ/d \] for a very general hypersurface $X\subset\PP^{n+1}$ of degree~$d$ and for $\frac n2<c<n$.
\end{corollary}

In particular, the integral Hodge conjecture fails for very general hypersurfaces $X\subset\PP^{n+1}$ of these degrees~$d$.

For $n=3$, the previous result from \cite{debarre-hulek-spandaw} (item \ref{it:6d} above) allows to disprove the integral Hodge conjecture for a set of degrees with density $\frac16$.
With our approach, we can actually show the failure of the integral Hodge conjecture for a set of degrees with density~$1$:
\begin{theorem}\label{thm:ihc}
Let $n\ge3$ be an integer.
Then there exists a set of degrees~$d$ with density $1$ such that the integral Hodge conjecture for very general hypersurfaces $X\subset\PP^{n+1}$ of degree~$d$ is false in every codimension~$c$ with $\frac n2<c<n$.
\end{theorem}

Theorems \ref{thm:main} and \ref{thm:ihc} as well as the known results \ref{it:ddd}, \ref{it:3dd}, and \ref{it:6d}
concern very general hypersurfaces. This might possibly exclude all hypersurfaces defined over number fields.
However, Totaro proved in \cite{totaro} that \ref{it:ddd}, \ref{it:3dd}, and \ref{it:6d} are in most cases also valid for
certain smooth hypersurfaces $X\subset\PP^4$ defined over $\QQ$.
Using Totaro's results, we show:
\begin{theorem}\label{thm:Q}
There exists a smooth hypersurface $X\subset\PP^4$ of degree $d\ge6$ defined over $\QQ$
such that the degree of every curve $C\subset X_{\overline\QQ}$ is divisible by $d$.
\end{theorem}

\subsection{General observations and notation}
Let $X\subset\PP^{n+1}$ be a very general hypersurface of degree~$d$
and let $\alpha=c_1(\Ocal_X(1))\in H^{1,1}(X,\ZZ)$ denote the hyperplane class.

By the Lefschetz hyperplane theorem, the restriction map
\[ H^i(\PP^{n+1},\ZZ)\to H^i(X,\ZZ) \]
is an isomorphism for $i<n$. Therefore,
\[ H^{2c}(X,\ZZ)=H^{c,c}(X,\ZZ)=\ZZ\cdot\alpha^c \quad\text{for $c<\tfrac n2$.} \]
Hence, Conjecture~\ref{conj:higher} is true for subvarieties $Z\subset X$ of codimension $c<\frac n2$.
Moreover, $Z^{2c}(X)$ is trivial for $c<\frac n2$.

For $\frac n2<c\le n$, Poincar\'e duality implies that
\[ H^{2c}(X,\ZZ)=H^{c,c}(X,\ZZ)=\ZZ\cdot\frac1d\alpha^c \quad\text{for $\tfrac n2<c\le n$.} \]
Therefore, no topological obstructions on the degree of subvarieties $Z\subset X$ of codimension~$c$ with $\frac n2<c\le n$ exist.
In this case, we can rephrase Conjecture~\ref{conj:higher} in terms of the group $Z^{2c}(X)$.
Since $\alpha^c$ is clearly algebraic, $Z^{2c}(X)$ is a quotient of $\ZZ/d$.
The order of $Z^{2c}(X)$ is given by the greatest common divisor of the degrees of all subvarieties $Z\subset X$ of codimension~$c$.
Hence, Conjecture~\ref{conj:higher} in codimension~$c$ for $\frac n2<c<n$ is equivalent to $Z^{2c}(X)\cong\ZZ/d$.
In particular, Theorem~\ref{thm:main} implies Corollary~\ref{cor:Z}.

Note that it suffices to prove Conjecture~\ref{conj:higher} for curves $C\subset X$ because every positive-dimensional subvariety $Z\subset X$ gives rise to a curve $C\subset X$ of the same degree after intersecting $Z$ with a suitable linear subspace.

If $X\subset\PP^{n+1}$ is very general, the order of the group $Z^{2n-2}(X)$ only depends on $n$ and $d$. We set $f_n(d):=|Z^{2n-2}(X)|$ for $n\ge3$ and $d\ge1$.
In other words,
\[ f_n(d) = \gcd\{\deg C\mid C\subset X \text{ curve} \} \;. \]
For $n=3$, this agrees with the definition of $f_3(d)$ given earlier.

We know that $f_n(d)\mid d$ for all degrees~$d$, and Conjecture~\ref{conj:higher} in degree~$d$ is equivalent to $f_n(d)=d$.

\subsection{Proof idea and overview}
Looking at the existing results towards Conjecture~\ref{conj:gh}, statement \ref{it:6d} seems to be more powerful than \ref{it:ddd} and \ref{it:3dd}.
However, the main idea for proving Theorem~\ref{thm:main} is to combine \ref{it:ddd}, \ref{it:3dd}, and \ref{it:6d} in order to show $d\mid f_3(d)$ for certain degrees~$d$.
This is based on the observation that $d\mid f_n(d_1)$ and $d\mid f_n(d_2)$ together imply $d\mid f_n(d_1+d_2)$, as can be seen by another degeneration argument.

In section~\ref{sect:trento}, we develop higher-dimensional analogues of the Trento examples \cite{trento-examples}.
These allow to carry out our approach in arbitrary dimension $n\ge3$.

In section~\ref{sect:main}, we will prove the following statement behind Theorem~\ref{thm:main}:
\begin{proposition}\label{prop:higher}
Let $n\ge3$ be an integer.
Then we have $f_n(d)=d$ if $d$ is coprime to $n!$ and the largest prime power $q$ dividing $d$ satisfies
\[ \textstyle \left(\binom n2-1\right)\cdot q^n+\left(n!-\binom n2\right)\cdot q^{n-1}+(2^n+1)\cdot n! \le d \;. \]
\end{proposition}
For $n=3$, the smallest degree~$d$ with this property is \[ d=5\cdot7\cdot11\cdot13=5005 \;. \]

In section~\ref{sect:asymptotic}, we will see that the positive integers~$d$ fulfilling the condition in Proposition~\ref{prop:higher} have positive density for all $n\ge3$, thus completing the proof of Theorem~\ref{thm:main}.

Finally, we prove Theorem~\ref{thm:ihc} in section~\ref{sect:ihc} and Theorem~\ref{thm:Q} in section~\ref{sect:Q}.

\section*{Acknowledgements}

I am grateful to my supervisor Stefan Schreieder for many helpful suggestions and comments concerning this work.
Further, I would like to thank the anonymous referee for useful remarks that improved the exposition.

This project received partial funding by the DFG project ``Topological properties of algebraic varieties'' (grant no.\ 416054549)
and by the ERC grant ``RationAlgic'' (grant no.\ 948066).
The final part of this project was carried out while the author was in residence at Institut Mittag-Leffler,
supported by the Swedish Research Council under grant no.\ 2016-06596.

\section{The Trento examples}\label{sect:trento}

In this section, we take a closer look at the Trento examples from \cite{trento-examples} by generalizing them to arbitrary dimension $n\ge3$.

All examples rely on the following lemma:
\begin{lemma}[Koll\'ar]\label{lem:kollar}
Let $n\ge3$ be an integer.
Suppose that there exists a smooth projective variety $Y$ of dimension $n$ with a very ample line bundle $L$ such that $L^n=d$ and $k\mid B\cdot L$ for every curve $B\subset Y$.
Then we have \[ k\mid n!\cdot f_n(d) \;. \]
\end{lemma}
\begin{proof}
We consider the embedding $Y\subset\PP^N$ given by the very ample line bundle $L$, and take a general linear projection
\[ \pi\colon Y\to\PP^{n+1} \;. \]
Then $\pi(Y)\subset\PP^{n+1}$ is a hypersurface of degree $L^n=d$.

By \cite{mather} (see also \cite{beheshti-eisenbud}), each fiber of $\pi$ has at most $n+1$ distinct points
(note that the fibers of $\pi$ might have a much larger degree than $n+1$ due to their non-reduced scheme structure, but for our argument we only need that they consist of at most $n+1$ points topologically).
Moreover, only finitely many fibers have exactly $n+1$ distinct points.
Hence, for every curve $C\subset\pi(Y)$, the curve $B=\pi^{-1}(C)^{\text{red}}\subset Y$ admits a finite surjective map $\pi|_B\colon B\to C$ of degree at most $n$.
Therefore, we have $B\cdot L\mid n!\cdot\deg C$ and thus $k\mid n!\cdot\deg C$.

Now if $X\subset\PP^{n+1}$ is a very general hypersurface of degree~$d$, every curve on $X$ specializes to a curve $C\subset\pi(Y)$ of the same degree.
For more details on this degeneration argument, see \cite[section~2]{soule-voisin}.
Since $k\mid n!\cdot\deg C$, it follows that $k\mid n!\cdot f_n(d)$.
\end{proof}
\begin{corollary}[Koll\'ar]\label{cor:ddd}
For all $n\ge3$ and $d\ge1$, we have \[ d\mid n!\cdot f_n(d^n) \;. \]
In particular, we have $d\mid f_n(d^n)$ if $d$ is coprime to $n!$.
\end{corollary}
\begin{proof}
We apply Lemma~\ref{lem:kollar} to $Y=\PP^n$ and $L=\Ocal_{\PP^n}(d)$.
\end{proof}
\begin{corollary}[Koll\'ar]\label{cor:3dd}
For all $n\ge3$ and $d\ge4$, we have \[ \textstyle d\mid n!\cdot f_n\left(\binom n2 d^{n-1}\right) \;. \]
In particular, we have $d\mid f_n\left(\binom n2 d^{n-1}\right)$ if $d$ is coprime to $n!$.
\end{corollary}
\begin{proof}
We take $Y=S\times\PP^{n-2}$ where $S\subset\PP^3$ is a very general surface of degree $d\ge4$.
On $Y$, we consider the very ample line bundle
\[ L=\pr_1^*\Ocal_S(1)\otimes\pr_2^*\Ocal_{\PP^{n-2}}(d) \;. \]
Then we have $L^n=\binom n2d^{n-1}$ (note that the factor $\binom n2$ is accidentally missing in \cite{trento-examples}).
If $B\subset Y$ is a curve, we obtain
\[ B\cdot L = (\pr_1)_*B\cdot\Ocal_S(1) + (\pr_2)_*B\cdot\Ocal_{\PP^{n-2}}(d) \equiv 0 \pmod d \;, \]
because $d$ divides the degree of the curve $\pr_1(B)\subset S$ by the Noether--Lefschetz theorem \cite{lefschetz}.
Therefore, Lemma~\ref{lem:kollar} implies the result.
\end{proof}
\begin{corollary}[van Geemen, Debarre--Hulek--Spandaw]\label{cor:6d}
For all $n\ge3$ and $d\ge2^n+1$, we have \[ d\mid(n-1)!\cdot f_n(n!\cdot d) \;. \]
In particular, we have $d\mid f_n(n!\cdot d)$ if $d$ is coprime to $(n-1)!$.
\end{corollary}
\begin{proof}
Let $(Y,L)$ be a very general polarized Abelian variety of dimension $n$ and type $(1,\ldots,1,d)$.
Then we have $L^n=n!\cdot d$.
It was shown in \cite{debarre-hulek-spandaw} that the line bundle $L$ is very ample.

Since $L$ is of type $(1,\ldots,1,d)$, we have
\[ c_1(L) = \dd x_1\wedge\dd x_2 + \cdots + \dd x_{2n-3}\wedge\dd x_{2n-2} + d\cdot\dd x_{2n-1}\wedge\dd x_{2n} \in H^2(Y,\ZZ) \]
for a suitable basis $\dd x_1,\ldots,\dd x_{2n}$ of $H^1(Y,\ZZ)$.
From this we see that
\[ \frac{c_1(L)^{n-1}}{(n-1)!} = \dd x_1\wedge\dd x_2\wedge\cdots\wedge\dd x_{2n-3}\wedge\dd x_{2n-2} + d\cdot\ldots \in H^{2n-2}(Y,\ZZ) \]
is not divisible by any integer larger than $1$.
If $Y$ is very general, the algebraic classes in $H^{2n-2}(Y,\ZZ)$ are rational multiples of $c_1(L)^{n-1}$,
and thus integral multiples of $c_1(L)^{n-1}/(n-1)!$.
Therefore, the degree of every curve $B\subset Y$ is divisible by $L^n/(n-1)!=nd$. Hence, Lemma~\ref{lem:kollar} gives $nd\mid n!\cdot f_n(n!\cdot d)$.
\end{proof}

\section{Proof of Proposition~\ref{prop:higher}}\label{sect:main}

We combine Corollaries \ref{cor:ddd}, \ref{cor:3dd}, and \ref{cor:6d} via the following simple observation:
\begin{lemma}\label{lem:deg}
If $d\mid f_n(d_1)$ and $d\mid f_n(d_2)$, then $d\mid f_n(d_1+d_2)$.
\end{lemma}
\begin{proof}
Let $C\subset X$ be a curve on a very general hypersurface $X\subset\PP^{n+1}$ of degree $d_1+d_2$.
By the same degeneration argument which we used in the proof of Lemma~\ref{lem:kollar},
\emph{every} hypersurface $X\subset\PP^{n+1}$ of degree $d_1+d_2$ contains a curve of the same degree as $C$.
In particular, we can choose $X=X_1\cup X_2$ to be the union of very general hypersurfaces $X_1,X_2\subset\PP^{n+1}$ of degrees $d_1$ and $d_2$, respectively.
Then every irreducible component of a curve $C\subset X_1\cup X_2$ lies on $X_1$ or $X_2$.
By assumption, the degrees of these components are divisible by $d$.
We conclude that $d\mid\deg C$.
\end{proof}

Now we can prove Proposition~\ref{prop:higher}.

\begin{proof}[Proof of Proposition~\ref{prop:higher}]
Since $d$ is a product of pairwise coprime powers of primes,
it suffices to show $q\mid f_n(d)$ for every prime power $q\mid d$.
By assumption, we have
\begin{equation}
\textstyle \left(\binom n2-1\right)\cdot q^n+\left(n!-\binom n2\right)\cdot q^{n-1}+(2^n+1)\cdot n! \le d \;. \label{eq:cond}
\end{equation}

We choose $i\in\left\{0,\ldots,\binom n2-1\right\}$ such that \[ d\equiv i\cdot q^n\pmod{\tbinom n2} \;. \]
This is possible because $\binom n2$ divides $n!$ and $q$ is coprime to $n!$.

Then we choose $j\in\{0,\ldots,n!-1\}$ such that \[ d\equiv i\cdot q^n+j\cdot q^{n-1}\pmod{n!} \;. \]
Our choice of $i$ implies that $j$ is divisible by $\binom n2$, so we have $j\le n!-\binom n2$.

By our choice of $j$, there exists an integer $k$ such that
\[ d = i\cdot q^n + j\cdot q^{n-1} + k\cdot n! \;. \]
Since $q\mid d$, we have $q\mid k$.
And from \eqref{eq:cond} we get $k\ge2^n+1$.

Now we have:
\begin{itemize}
\item $q\mid f_n(q^n)$ by Corollary~\ref{cor:ddd}
\item $q\mid f_n(\binom n2q^{n-1})$ by Corollary~\ref{cor:3dd}
\item $k\mid(n-1)!\cdot f_n(k\cdot n!)$ by Corollary~\ref{cor:6d} and thus $q\mid f_n(k\cdot n!)$
\end{itemize}
Combining these results via repeated usage of Lemma~\ref{lem:deg}, we obtain
\[ q\mid f_n\left(i\cdot q^n+\frac j{\binom n2}\cdot\tbinom n2q^{n-1}+k\cdot n!\right) = f_n(d) \;. \qedhere \]
\end{proof}

\begin{remark}
Using only Corollaries \ref{cor:ddd} and \ref{cor:6d},
one can show a weaker statement where \eqref{eq:cond} is replaced by the assumption
\[ (n!-1)\cdot q^n+(2^n+1)\cdot n!\le d \;. \]
It turns out that this stronger condition results in a set of degrees $d$ with the same density.
Therefore, Corollary~\ref{cor:3dd} is not strictly necessary to obtain Theorem~\ref{thm:main}.
However, only together with Corollary~\ref{cor:3dd} we can prove Conjecture~\ref{conj:gh} for $d=5005$.
\end{remark}

\section{Some analytic number theory}\label{sect:asymptotic}

A set $A$ of positive integers has \emph{positive density} if
\[ \liminf_{m\to\infty} \frac{|A\cap\{1,\ldots,m\}|}m > 0 \;. \]

In this section, we want to prove the following:
\begin{proposition}\label{prop:asymptotic}
Let $n\ge1$ be an integer and $\lambda>0$ a real number.
Then the positive integers~$d$ coprime to $n!$ such that the largest prime power dividing $d$ is not larger than $\lambda\cdot d^{1/n}$ have positive density.
\end{proposition}
Together with Proposition~\ref{prop:higher}, this will complete the proof of Theorem~\ref{thm:main},
since for $d\gg0$ the condition $q\le\binom n2^{-1/n}d^{1/n}$ implies
\[ \textstyle \left(\binom n2-1\right)\cdot q^n+\left(n!-\binom n2\right)\cdot q^{n-1}+(2^n+1)\cdot n! \le d \;, \]
so Proposition~\ref{prop:higher} applies to a set of degrees~$d$ with positive density.

We use the following easy lemma on the distribution of prime powers:
\begin{lemma}\label{lem:Pi}
Let $\Pi(m)$ denote the number of prime powers $\le m$. Then
\[ \frac{\Pi(m)}m \overset{m\to\infty}\longrightarrow 0 \;. \]
\end{lemma}
\begin{proof}
By the prime number theorem, we have
\[ \frac{\pi(m)}m \overset{m\to\infty}\longrightarrow 0 \;, \]
where $\pi(m)$ counts the prime numbers $\le m$.
Now if $p^e\le m$ is a prime power with $e\ge2$,
we have $e\le\log_2 m$ and $p\le\sqrt m$, so we conclude by noting that
\[ \frac{\log_2 m\cdot\sqrt m}m \overset{m\to\infty}\longrightarrow 0 \;. \qedhere \]
\end{proof}
We also need the following consequence of Mertens' theorem:
\begin{lemma}\label{lem:mertens}
We have
\[ \sum_{x^{1/n}<p\le x}\frac 1p \overset{x\to\infty}\longrightarrow\log n \;, \]
where the sum runs only over prime numbers $p$.
\end{lemma}
\begin{proof}
By \cite{mertens}, there exists a constant $C$ such that
\[ \sum_{p\le x}\frac1p - \log\log x \overset{x\to\infty}\longrightarrow C \;. \]
Since $\log\log x-\log\log x^{1/n}=\log n$, we conclude that
\[ \sum_{x^{1/n}<p\le x}\frac1p = \sum_{p\le x}\frac1p-\sum_{p\le x^{1/n}}\frac1p \overset{x\to\infty}\longrightarrow \log n \;. \qedhere \]
\end{proof}
Now we can prove Proposition~\ref{prop:asymptotic}.
\begin{proof}[Proof of Proposition~\ref{prop:asymptotic}]
Let $\alpha>1$ be a real number.
Dickman proved in \cite{dickman} that the positive integers~$d$ whose largest prime divisor is not larger than $d^{1/\alpha}$ have density $\rho(\alpha)$, where $\rho$ denotes the Dickman function.
More generally, this result was proven for arithmetic progressions in \cite{buchstab}.
Therefore, the positive integers~$d$ coprime to $n!$ whose largest prime divisor is not larger than $d^{1/\alpha}$ have density \[ \frac{\varphi(n!)}{n!}\cdot\rho(\alpha) \;, \] where $\varphi$ denotes Euler's totient function.
Since $\rho$ is continuous, it follows that the positive integers~$d$ coprime to $n!$ whose largest prime divisor is not larger than $\lambda\cdot d^{1/n}$ have density $\frac{\varphi(n!)}{n!}\cdot \rho(n)>0$.

In other words, Proposition~\ref{prop:asymptotic} holds if we replace `prime power' by `prime number'.
Hence, it suffices to show that the positive integers~$d$ divisible by a prime power $q=p^e>\lambda\cdot d^{1/n}$ with $e\ge2$ have density $0$.

For a given $x$, let us consider the number $N(x)$ of positive integers $d\le x$ with this property.
Any such $d$ can be written as \[ d=q\cdot r \;, \] where $q=p^e\ge\lambda\cdot d^{1/n}$ is a prime power with $e\ge2$.
For fixed $q\le x^{1/n}$, there are at most $\lambda^{-n}q^{n-1}$ possibilities for $d$ because \[ r=\frac dq\le\frac{\lambda^{-n}q^n}q=\lambda^{-n}q^{n-1} \;. \]
For fixed $q>x^{1/n}$, there are at most $\frac xq$ possibilities for $d$ because \[ r=\frac dq\le\frac xq \;. \]
Together we obtain the upper bound
\[ N(x)\ \le\ \lambda^{-n}\cdot\sum_{q\le x^{1/n}}q^{n-1}\ +\ x\cdot \sum_{x^{1/n}<q\le x}\frac1q \;, \]
where both sums run only over prime powers $q=p^e$ with $e\ge2$.

Using Lemma~\ref{lem:Pi}, we get
\[ \frac1x\cdot\sum_{q\le x^{1/n}}q^{n-1}\ \le\ \frac1x\cdot\Pi(x^{1/n})\cdot\left(x^{1/n}\right)^{n-1}=\frac{\Pi(x^{1/n})}{x^{1/n}} \overset{x\to\infty}\longrightarrow 0 \;, \]
so in order to prove $\frac{N(x)}x\to0$ for $x\to\infty$, it remains to show that
\[ \sum_{x^{1/n}<q\le x}\frac1q \overset{x\to\infty}\longrightarrow 0  \;. \]
For $q=p^e\le x$, we have $e\le\log_2 x$ and thus
\begin{align*}
\limsup_{x\to\infty}\sum_{x^{1/n}<q\le x}\frac 1q
&= \limsup_{x\to\infty}\sum_{e=2}^{\lfloor\log_2x\rfloor}\sum_{x^{1/n}<p^e\le x}\frac1{p^e} \\
&\le \limsup_{x\to\infty}\left(\sum_{e=2}^n \sum_{x^{1/en}<p\le x^{1/e}}\frac{x^{-\frac{e-1}{en}}}p + \sum_{e=n+1}^{\lfloor\log_2x\rfloor}\frac{x^{1/e}}{x^{1/n}} \right) \;. \\
\intertext{Here we used $\frac1{p^e}=\frac{p^{-(e-1)}}p\le\frac{x^{-\frac{e-1}{en}}}p$ for $p>x^{1/en}$ if $2\le e\le n$, and $\frac1{p^e}\le\frac1{x^{1/n}}$ for $p>x^{1/en}$ if $e\ge n+1$.
Applying Lemma~\ref{lem:mertens} for each $2\le e\le n$, and using $x^{1/e}\le x^{1/(n+1)}$ for $e\ge n+1$, we obtain}
\cdots &\le \limsup_{x\to\infty}\left(\sum_{e=2}^n\frac{\log n}{x^{\frac{e-1}{en}}}+\frac{\log_2x}{x^{\frac1{n(n+1)}}}\right) = 0 \;. \qedhere
\end{align*}
\end{proof}
\begin{remark}
The proof shows that the density in Theorem~\ref{thm:main} amounts to \[ \frac{\varphi(n!)}{n!}\cdot\rho(n) \;. \]
For example, the density for $n=3$ is $\frac13\cdot\rho(3) \approx 1.6\%$.
\end{remark}

\section{Failure of the integral Hodge conjecture}\label{sect:ihc}

By the work of Koll\'ar \cite{trento-examples}, hypersurfaces provide an example for varieties where the integral Hodge conjecture fails due to a non-torsion cohomology class.
Theorem~\ref{thm:ihc} says that this counterexample works for almost all degrees~$d$ (in the sense of density).
\begin{proof}[Proof of Theorem~\ref{thm:ihc}]
The failure of the integral Hodge conjecture in degree~$d$ is equivalent to $f_n(d)\ne1$.
Hence, we need to show that the positive integers~$d$ with $f_n(d)\ne1$ have density $1$.
If $d$ has a prime divisor $p$ coprime to $n!$ such that
\[ \textstyle \left(\binom n2-1\right)\cdot p^n + \left(n!-\binom n2\right)\cdot p^{n-1} + (2^n+1)\cdot n!\le d \;, \]
then the proof of Proposition~\ref{prop:higher} shows that $p\mid f_n(d)$.
Therefore, for every prime $p>n$ all sufficiently large multiples~$d$ of $p$ satisfy $f_n(d)\ne1$.

For a given $\varepsilon>0$, we can find distinct primes $p_1,\ldots,p_N>n$ such that
\[ \frac1{p_1}+\cdots+\frac1{p_N} > \frac1\varepsilon \]
since the sum of the reciprocals of all primes diverges.
We know from the previous paragraph that for $d\gg0$, we might have $f_n(d)=1$ only if $d$ is not divisible by any of the primes $p_1,\ldots,p_N$.
Hence, the density of these $d$ is at most
\[ \left(1-\frac1{p_1}\right)\cdots\left(1-\frac1{p_N}\right)<\frac1{\left(1+\frac1{p_1}\right)\cdots\left(1+\frac1{p_N}\right)}<\frac1{\frac1{p_1}+\cdots+\frac1{p_N}}<\varepsilon \;. \]
This concludes the proof.
\end{proof}
\begin{remark}
For any degree $d\ge1$, there exist special smooth hypersurfaces $X\subset\PP^{n+1}$ of degree~$d$ which do satisfy the integral Hodge conjecture in every codimension~$c$ with $\frac n2<c<n$. For example, we can take the Fermat hypersurface \[ \{x_0^d+\cdots+x_{n+1}^d=0\}\subset\PP^{n+1} \;, \]
since it contains an $(n-c)$-dimensional linear subspace for any $\frac n2<c<n$.
\end{remark}
\begin{remark}
There are infinitely many degrees~$d$ for which we are not able to disprove the integral Hodge conjecture.
In particular, this problem remains open when $d$ is a prime number,
in which case the failure of the integral Hodge conjecture is equivalent to Conjecture~\ref{conj:higher}.
\end{remark}

\section{Example over $\QQ$}\label{sect:Q}

The basic idea in \cite{totaro} is to replace the original degeneration arguments with degenerations to positive characteristic.
To prove Theorem~\ref{thm:Q}, we apply this idea to the proof of Proposition~\ref{prop:higher} and use some of Totaro's results.
\begin{proposition}
There exists a smooth hypersurface $X\subset\PP^4$ of degree $d=7\cdot13\cdot19\cdot31=53599$ defined over $\QQ$
such that the degree of every curve $C\subset X_{\overline\QQ}$ is divisible by $d$.
\end{proposition}
\begin{proof}
We will show the following lemma:
\begin{lemma}\label{lem:Fp}
Let $q$ be any of the four prime divisors of $d$.
Then there exists a prime~$p$ and a hypersurface $Y\subset\PP^4_{\FF_p}$ of degree~$d$ such that the degree of every curve $C\subset Y_{\overline{\FF_p}}$ is divisible by $q$.
Moreover, $p$ can be chosen such that finitely many given primes are avoided.
\end{lemma}
Once this lemma is proven, we proceed as follows:
We ensure that the primes~$p$ for each prime divisor $q\mid d$ are pairwise different.
Then we use the Chinese remainder theorem to construct a smooth hypersurface $X\subset\PP^4$ defined over $\QQ$ that simultaneously specializes to all four hypersurfaces $Y\subset\PP^4_{\FF_p}$ from Lemma~\ref{lem:Fp}.
This $X$ satisfies our claim.
\end{proof}
\begin{proof}[Proof of Lemma~\ref{lem:Fp}]
Since $q\equiv1\pmod6$, we can write $d=q^3+6k$ for some integer $k$. Note that $k\ge38$.
As in the proof of Lemma~\ref{lem:deg}, we want to take $Y=Y_1\cup Y_2$, where $Y_1,Y_2\subset\PP^4_{\FF_p}$ are two hypersurfaces of degrees $q^3$ and $6k$, respectively, such that every curve $C\subset (Y_i)_{\overline{\FF_p}}$ has degree divisible by $q$.

We first construct $Y_1$. A priori, \cite[Corollary~4.2]{totaro} only gives hypersurfaces $Y_1\subset\PP^4$ over $\overline{\FF_p}$ with this property for every $p>q^3$.
However, as in the proofs of \cite[Lemma~5.1]{totaro} and \cite[Theorem~6.1]{totaro}, we can apply \cite[Lemma~4.3]{totaro} to $\PP^3_\QQ$ (polarized by $\Ocal_{\PP^3}(q)$) to get a rational map to $\PP^4_\ZZ$ and obtain hypersurfaces $Y_1\subset\PP^4$ over $\FF_p$ after excluding finitely many primes~$p$.

Now we construct $Y_2$. The proof of \cite[Theorem~6.1]{totaro} yields a prime~$p$ and a hypersurface $Y_2\subset\PP^4_{\FF_p}$ of degree $6k$ such that $k\mid 6\cdot\deg C$ for every curve $C\subset(Y_2)_{\overline{\FF_p}}$. Since $k$ is a multiple of $q$ and $q$ is coprime to $6$, it follows that $q\mid\deg C$. Furthermore, we can guarantee that $p$ is different from finitely many given primes (including also the primes where the construction of $Y_1$ does not work) by doing the argument of \cite[Theorem~6.1]{totaro} over $\ZZ[1/P]$ instead of $\ZZ$, where $P$ is the product of these finitely many primes.
\end{proof}

\begin{remark}
For simplicity, we gave only one specific example over $\QQ$.
The above argument obviously works for other values than $d=53599$ as well.
\end{remark}

\bibliographystyle{alphaurl}
\bibliography{gh}

\end{document}